\documentclass[12pt,a4paper]{article}
      \usepackage[english]{babel}
      \usepackage{setspace}
	  \usepackage{amsmath}
      \usepackage{amstext}
      \usepackage{amssymb}
      \usepackage{amsthm}
      \usepackage{amsfonts}
      \usepackage[T1]{fontenc}
      \usepackage{latexsym}
      \usepackage{psfrag}
	  \newtheorem*{thm*}{Theorem}
      \newtheorem{thm}{Theorem}[section]

      \newcommand{\N}{\mathbb N}
      \newcommand{\R}{\mathbb R}

      \newcommand{\ds}{\displaystyle}
      \newcommand{\U}{\mathbb U}

      \newenvironment{appli}{\left( \begin{array}{ccc}}{\end{array} \right)}

      \newcommand{\ba}{\begin{appli}}
      \newcommand{\ea}{\end{appli}}
\begin{document}
%\doublespace

\title{Computing the complexity of the relation of isometry between separable Banach spaces}
\author{Julien Melleray}
%\thanks{\textit{2000 Mathematical Subject Classification}: 03E15.} \\
%\textit{Keywords}: Urysohn universal space, isometry, definable equivalence relation}
\date{}
\maketitle
\begin{abstract} \noindent We compute here
the Borel complexity of the relation of isometry between separable Banach spaces, using results of 
Gao-Kechris \cite{gaokec} and Weaver \cite{weaver}.
\end{abstract}

\begin{section}{Introduction}
Over the past fifteen years or so, the theory of complexity of Borel equivalence relations has been a very active field of 
research; in this paper, we compute the complexity of a relation of geometric nature, the 
relation of (linear) isometry between separable Banach spaces. 
Before stating precisely our result, we begin by recalling  the basic facts and definitions that we need 
in the following of the article; we refer the reader to \cite{Kechris1} 
for a thorough introduction to the concepts and methods of descriptive set theory.\\

{\bf (A) Notations, definitions.}\\
In this article, the letters $X$ and $Y$ always refer to separable metric spaces.\\
We say that a metric space $(X,d)$ is \textit{Polish} if $(X,d)$ is separable and complete; we often 
forget $d$ and write it simply $X$. \\
We call \textit{standard Borel space} any pair $(X,\Sigma)$, where $X$ is a set and $\Sigma$ is 
a $\sigma$-algebra on $X$ which is isomorphic to the $\sigma$-algebra formed by the Borel subsets of $[0,1]$. Any 
uncountable Polish space $X$, endowed with the $\sigma$-algebra whose elements are the Borel subsets of $X$, is a standard Borel space. \\
If $X$ and $Y$ are standard Borel spaces, a map $f \colon (X,\Sigma_X) \to (Y,\Sigma_Y)$ is said to be \textit{Borel} if 
$f^{-1}(A) \in \Sigma_X$ for all $A \in \Sigma_Y$. Borel maps are closed under composition. \\
If $X$ is a Polish metric space, we let ${\mathcal F}(X) =\{F \subset X \colon F \mbox{ is closed }\}$.\\
We endow it with the $\sigma$-algebra generated  by the sets $\{F \in {\mathcal F} (X) \colon F \cap U \neq \emptyset\}$, 
where $U$ varies over the open subsets of $X$. This is called the \textit{Effros Borel structure}.\\
If $X$ is uncountable then ${\mathcal F}(X)$, equipped with the Effros Borel structure, is a standard 
Borel space.\\
If $X$ is a standard Borel space and $A \subset X$, we say that $A$ is \textit{analytic} if 
there exists a standard Borel space $Y$ and a Borel map $f \colon Y \to X$ such that $A=f(Y)$.\\
We say that an equivalence relation $E$ defined on the standard Borel space $X$ is Borel (resp. analytic) if 
it is a Borel (resp. analytic) subset of $X \times X$, endowed with its natural Borel structure.\\
To avoid confusions, let us specify that an \textit{isometry} is a bijective map $f \colon X \to Y$ such that 
$d(f(x),f(x'))=d(x,x')$; if $f \colon X \to Y$ is distance-preserving but not onto, 
then we say that $f$ is an \textit{isometric embedding}.\\ 
If $X$ is a Polish metric space, we let $Iso(X)$ denote its isometry group.\\

{\bf (B) Classification problems}\\
In our setting, a classification problem is, given a "definable" (Borel, analytic...) equivalence relation 
$E$ on a standard Borel space $X$, 
to find invariants for this relation, i.e a set $I$ and a map $f \colon X \to I$ such that 
$(xEy) \Leftrightarrow (f(x)=f(y))$.\\
Of course, for this to be of interest, both $I$ and $f$ have to be as concrete as possible; 
we refer the reader to \cite{gaokec} or \cite{KEchris2} for a detailed introduction about classification problems, 
and bibliographical references on this subject.\\
Since some classification problems have already been well-studied, it is natural to look for a way to compare the 
complexities of classification problems:
intuitively, an equivalence relation $E$ is simpler than another equivalence relation $E'$ when knowing a classification 
for $E'$ is enough to obtain a classification for $E$.\\
Formally, one says that a relation $E$, defined on a standard Borel space $X$, \textit{Borel reduces to}
a relation $E'$ on a standard Borel space $X'$ if there exists a Borel map $f \colon X \to X'$
such that
$$(xEy) \Leftrightarrow (f(x) E' f(y)) \ .$$
This indeed implies that "$E$ is simpler than $E'$", in the sense explained above: composing by $f$, any 
classification of $E'$ is enough to define a classification of $E$ (asking that $f$ be Borel means that it is 
also not too much harder to compute the invariants for $E$ when knowing those for $E'$). \\
We note $E \leq_B E'$ if $E$ Borel reduces to $E'$, and $E \sim_B E'$ if $E$ and $E'$ Borel reduce to one another
(in other worlds, the classification problems for $E$ and $E'$ have the same complexity).\\
The relation $\leq_B$ induces a hierarchy on the levels of complexity of equivalence relations; in this article, 
we concern ourselves with the relation of isometry between Banach spaces, and compute its place in 
this hierarchy. To explain how this relation fits in our frame, and to state precisely our main result, we need to recall 
some properties a a remarkable Polish space, the \textit{universal space of Urysohn}.\\

{\bf (C) Urysohn's universal metric space.}\\
Up to isometry, Urysohn's universal metric space $\U$, first constructed by Urysohn in \cite{Urysohn}, 
is the only Polish metric space with has the following two properties: \\
- $\U$ is universal, which means that any Polish metric space is isometric to a (necessarily closed) subset of $\U$;\\
- $\U$ is $\omega$-homogeneous, i.e any isometry between two finite subsets $F_1,F_2$ of $\U$ extends to an isometry of $\U$. \\
For more informations about this space, and bibliographical references, we refer the reader to \cite{gaokec} or \cite{melleray3}. \\

Here, we use the Urysohn space because of results by Gao and Kechris \cite{gaokec}. 
To state these, we first need to point out that, since any Polish 
metric space is isometric to some closed set $F \in {\mathcal F}(\U)$, one 
may consider ${\mathcal F}(\U)$ (with the Effros Borel structure), as being the (Borel) space of Polish spaces. \\
One easily checks that the relation of isometry $\backsimeq_i$ (defined on ${\mathcal F}(\U)$, endowed with the 
Effros Borel structure) is analytic, where 
$$ (P \backsimeq_i P') \Leftrightarrow (P \mbox{ and } P' \mbox{ are isometric})\  .$$ 
To compute the exact complexity 
of this relation, Gao and Kechris considered the relation $\backsimeq_i^{\U}$ defined, for $P,P' \in {\mathcal F}(\U)$, by 
$$(P \backsimeq_i^{\U} P') \Leftrightarrow (\exists \varphi \in Iso(\U) \ \varphi(P)=P')\ . $$
Using a variation of Kat\v{e}tov's construction of $\U$ (cf. \cite{Katetov}), they proved that 
$(\backsimeq_i) \sim_B (\backsimeq_i^{\U})$, and that $\backsimeq_i^{\U}$ is Borel bireducible to 
the universal relation for relations induced by a Borel action of a Polish group. \\
We use this result to compute the complexity of the relation of (linear) isometry between separable Banach spaces.\\

{\bf (D) Organization of the article and statement of the main theorem.} \\
We begin by briefly presenting the basic facts of the theory of Lipschitz-free Banach spaces, and remark that results of 
Godefroy and Kalton \cite{gilleskalton} about these spaces are enough to show that the uniquely determined closed linear 
span of $\U$ (cf. \cite{holmes}) is a universal Banach space up to linear isometry, thus answering a question of Holmes \cite{holmes}. \\
This Banach space is actually isometric to $F(\U)$, the Lipschitz-free space over $\U$ (see the next section for definitions). \\
If $B$ is a Banach space, $\{F \in {\mathcal F}(B) \colon F \mbox{ is a subspace}\}$ is a Borel subset of 
${\mathcal F}(B)$; we may thus consider $\{F \in {\mathcal F}(F(\U)) \colon F \mbox{ is a subspace}\}$ as the (Borel) space 
of all separable Banach spaces. 
Then one again checks that the relation of linear isometry between separable Banach spaces is analytic; 
since any isometry between Banach spaces is affine, it is clear that this relation Borel reduces to $\backsimeq_i$. \\
We use results of Weaver \cite{weaver} to prove that the relation of isometry between Polish metric 
spaces Borel reduces to the relation of isometry between separable Banach spaces; 
this, added to the results of Gao and Kechris described above, 
enables us to compute the Borel complexity of the relation of isometry between separable Banach spaces.

\begin{thm*}
The relation of (linear) isometry between separable Banach spaces is Borel bireducible with the universal relation 
for relations induced by a Borel action of a Polish group.
\end{thm*}

\emph{Acknowledgements.} If not for conversations with several people, I would not have heard about many 
of the results used in this article, nor would I have understood them. 
I would especially like to thank Lionel Nguyen Van The, who told me about the article of Holmes \cite{holmes}; 
Valentin Ferenczi, who explained to me the links between this article and the theory of Lipschitz-free Banach spaces, and told me about 
the Godefroy-Kalton theorem; and Gilles Godefroy, who kindly took time to discuss with me the theory of Lipschitz-free 
Banach spaces. I am very grateful to all of them.
\end{section}

\begin{section}{Lipschitz-free Banach spaces }
In this section we briefly detail the basic definitions and properties of Lipschitz-free Banach spaces; 
we follow the second chapter of \cite{weaver}. The interested reader may find more informations about these 
spaces in \cite{weaver} and \cite{gilleskalton}. \\

If $(X,d,e)$ is a pointed metric space, one lets
$Lip_0(X,d,e)$ denote the space of Lipschitz functions on $X$ that map $e$ to $0$. \\

One defines a norm on $Lip_0(X,d,e)$ by setting
$$||f||= \inf\{k\in \R \colon f \mbox{ is } k-\mbox{Lipschitz} \} \ \ .$$
It is worth noting that, if one chooses another basepoint $e'$, then
$Lip_0(X,d,e)$ and $Lip_0(X,d,e')$ are isometric, one possible isometry being given by the mapping 
$f \mapsto f - f(e')$. \\
In the following, when no confusion is possible, we forget about $d$ and $e$
and simply write $Lip_0(X)$; we write $[Lip_0(X)]_1$ to denote the closed unit ball of $Lip_0(X)$. \\

If $(X,d)$ is a metric space, we say that $m \colon X \to
\R$ is a \textit{molecule} if $m$ has a finite support, and $\sum_{x \in X} m(x)=0$. \\
For $p,q \in X$, one may define a molecule $m_{pq}$ by setting
$\ds{m_{pq}=\chi_{\{p\}}-\chi_{\{q\}}}$, where $\ds{\chi_X}$ stands for the characteristic function of $X$. \\
For any molecule $m$, one may find points $p_i,q_i \in X$ and reals $a_i$ such that 
$ \ds{ m= \sum_{i=1}^n a_i m_{p_iq_i}}$. \\
 We let $\ds{ ||m ||= \inf \{ \sum_{i=1}^n |a_i| d(p_i,q_i)\colon m= \sum_{i=1}^n a_i m_{p_iq_i} \}\ .}$ 

Then $||.||$ is a seminorm on the space of molecules; we call \textit{Lipschitz-free space over} $X$, and note
$F(X)$, the completion (relative to this seminorm) of the space of molecules modulo null vectors 
(there are actually no null vectors, as we will see shortly). 
This space is also known in the litterature as the \emph{Arens-Eells space} of $X$. \\
The following fact is the basis of the theory of Lipschitz-free Banach spaces:\\

{\bf Fact.}(see e.g \cite{weaver}) The spaces $F(X)^*$ and $Lip_0(X)$ are isometric. \\

The natural isometry $T\colon F(X)^* \to Lip_0(X)$ is defined by\\
 $$  \forall x \in X \ (T\phi)(x)=\phi(m_{xe}) \ .$$
(here $e$ is any point in $X$; recall that, if $e \neq e' \in X$, then $Lip_0(X,d,e)$ and
$Lip_0(X,d,e')$ are isometric).\\
The inverse $S$ of $T$ is defined by the formula
$$ (Sf)(m)= \sum_{x \in X} f(x)m(x) \ .$$
Therefore, the Hahn-Banach theorem implies that, for any molecule $m$, one has
$$||m||=\sup \big \{\sum_{x\in X} f(x) m(x) \colon f \in [Lip_0(X)]_1 \big \} \ .$$
It is interesting to notice that this means that $||m||$ is determined by the distances between 
points in the support of $m$ and the values of $m$ on its support. Indeed, if $Y \subset X$ is 
a subspace of the metric space $X$, it is easy to see that any $1$-Lipschitz real-valued map on $Y$ extends 
to a $1$-Lipschitz real-valued map on $X$. Therefore, for any point $e \in X$, one has, for 
any $S$ that contains the support of $m$, that
$$||m||=\sup\big \{\sum_{x\in X} f(x) m(x) \colon f \in [Lip_0(S\cup\{e\},d,e)]_1 \big \} $$

This implies that $||.||$ is in fact a norm on the space of molecules. \\
Indeed, let $e \in X$, and $m$ be a non-zero molecule;  then one has \\
$m = \sum_{i=1}^n a_i m_{p_ie}$ for points $p_i \in X \setminus\{e\}$ 
(which we may assume to be pairwise distinct) and reals 
$a_i \neq 0$. \\
We may find some $\varepsilon >0$ such that $ \varepsilon < \min \{d(p_i,p_j) \colon i \neq j\}$. \\
Let now $f(e)=0$, $f(p_i)= \varepsilon \frac{|a_i|}{a_i}$; $f$ is $1$-Lipschitz.  \\
Then the formula above implies that
$$||m|| \geq \sum_{x \in X} f(x)m(x)= \sum_{i=1}^n \varepsilon |a_i| >0\ \ .$$
This proof also shows that the family $\{m_{xe}\}_{x \in X}$ is linearly independent.\\

Furthermore, one checks easily that $||m_{xy}||=d(x,y)$ for all $x,y \in X$.\\
Hence, if $e$ is any point in $X$, then the mapping $x \mapsto m_{xe}$ 
is an isometric embedding of $X$ in $F(X)$, such that the closed linear span of the image of $X$
is equal to the whole space $F(X)$.\\
In the following, we will be mainly interested in the Lipschitz-free space over the Urysohn space $\U$. \\

In \cite{holmes}, M.R Holmes, following earlier work of Sierpinski \cite{sierpinski} on isometric embeddings 
of $\U$ in Banach spaces, proved a very surprising result, which we state below.

\begin{thm*} (Holmes)
If $\U$ is isometrically embedded in a Banach space $B$, and $0 \in
\U$, then one has, for all $x_1,\ldots,x_n \in \U$ and $\lambda_1,\ldots,\lambda_n \in \R$:
$$
||\sum_{i=1}^n \lambda_i x_i|| = \sup \big \{ \big |\sum_{i=1}^n \lambda_i f(x_i)
\big | \colon f \in \big[Lip_0(\{x_1,\ldots,x_n \} \cup \{0\})\big]_1 \big
\} \ .
$$
\end{thm*}

This theorem has remarkable consequences: assume that $X,X'$ are isometric to $\U$, 
and that $0 \in X \subset B$, $0 \in X' \subset B'$, where $B$ and $B'$ are 
Banach spaces. \\
Then any isometry $\varphi \colon X \to X'$ mapping $0$ to $0$ extends to a linear isometry $\tilde{\varphi}$ 
which maps the closed linear span of $X$ (in $B$) onto the closed linear span of $X'$ (in $B'$): to see that, 
one simply has to check that the mapping 
$\tilde{\varphi} \colon \sum_{i=1}^n \lambda_i x_i \mapsto \sum_{i=1}^n \lambda_i \varphi(x_i)$ is an isometry from 
the linear span of $X$ to that of $X'$, and this is 
a direct consequence of the theorem of Holmes quoted above.\\
Also, if $P \subset X$ is a Polish metric space containing $0$, 
then Holmes' result shows that the closed linear span of $P$ (in $B$) is isometric to the Lipschitz-free space over $P$.\\
In particular, if $\U$ is embedded in a Banach space $B$ in such a way that $0 \in \U$, and the closed linear span of 
$\U$ is $B$, then $B$ is isometric to the Lipschitz-free Banach space over $\U$.\\

Holmes obviously did not know about the theory of Lipschitz-free Banach spaces; he called the unique (up to isometry) Banach 
space described above "the uniquely 
determined closed linear span of $\U$". After noticing that any separable metric space 
isometrically embeds into this space (since, as any Polish space, it isometrically embeds into $\U$!), 
he asked the following question (which is also the question n. 997 of \cite{openproblems}):\\

Is it true that any separable Banach space admits a \textit{linear} isometric embedding in the uniquely determined closed 
linear span of $\U$?\\

A theorem of Godefroy and Kalton \cite{gilleskalton} shows that the answer to this question is positive: 
indeed, they show that if $B,B'$ are separable Banach spaces and $B$ embeds isometrically in $B'$, then
$B$ embeds linearly isometrically in
$B'$. In particular, if $X$ is a universal Polish space up to isometry, then $F(X)$ 
is a universal separable Banach space up to 
linear isometry; in particular, $F(\U)$ is universal up to linear isometry. \\

{ \bf Remark: }If $X$ is a closed subset of the separable Banach space $B$, then the mapping 
$$F \in {\mathcal F}(X) \mapsto \overline {span (F)} \in {\mathcal F}(B)$$
is Borel (where both ${\mathcal F}(X)$ and ${\mathcal F}(B)$ are endowed with the Effros Borel structure).\\
So, if one identifies the class of Polish spaces to the set of closed 
subsets of $\U$ containing $0$, and the class of separable Banach spaces to the set 
of closed subspaces of $F(\U)$, then we may see the mapping $X \mapsto F(X)$ as a Borel mapping 
between two standard Borel spaces.\\

In our context, one question about Lipschitz-free Banach spaces is of special interest: 
if $X$ is a Polish metric space, how much of its metric structure is "encoded" in $F(X)$? In other words,
if one knows that $X,Y$ are Polish metric spaces such that $F(X),F(Y)$ are isometric, 
can we find a relation between the metric structures of $X,Y$ ?\\

Before saying more about this, we need a definition:\\
We say that $f \colon X \to Y$ is a \textit{dilatation} if there exists $\lambda >0$ such that 
$$d(f(x),f(x'))= \lambda d(x,x') \ .$$ 
It is easy to see that if there is a dilatation from $X$ onto $Y$ then $F(X)$ and $F(Y)$ 
are isometric. The converse is false in general, but beautiful results of Weaver imply that it holds 
for a rather large class of spaces; 
what is interesting for us is that the relation of isometry between all Polish spaces reduces to that 
of isometry between spaces in the aforementioned subclass.\\

Weaver \cite{weaver} says that a Polish space $P$ is \textit{concave} if,
for all $p \neq q \in P$, the molecule $\frac{m_{pq}}{d(p,q)}$ is an extreme point in the unit ball of $ F(X)^{**}$ 
(here we use the canonical embedding of $F(X)$ into its bidual).\\

Then Weaver proves the following result :
\begin{thm*} \emph{(}Weaver \cite{weaver}\emph{)} Let $P$ and $P'$ be two concave Polish metric spaces, and assume that $F(P)$ and
$F(P')$ are isometric. Then there exists a dilatation from $P$ onto
$P'$. 
\end{thm*}

He also shows that the class of concave Polish metric spaces is fairly large:
\begin{thm*} \emph{(}Weaver \cite{weaver}\emph{)} Let $(P,d)$ be a Polish metric space. Then $(P,\sqrt d)$ is concave.
\end{thm*}
(He actually proves that $(P,d^{\alpha})$ is concave for any $\alpha \in ]0,1[$; 
we will only need this fact for $\alpha=\frac{1}{2}$) 

Intuitively, by replacing $d$ by $\sqrt{d}$ (which is easily checked to be a complete distance, compatible 
with the topology of $P$), one has "uniformly eliminated" the equality case in the triangle inequality, and 
this fact is enough to study precisely the structure of the isometries of $F(P)$.\\

A very simple, yet very important for our constructions, fact is that, if $(P,d)$ and $(P',d')$ are 
two metric spaces, then $(P,d)$ and $(P',d')$ are isometric if, and only if, $(P,\sqrt{d})$
and $(P',\sqrt{d'})$ are isometric.\\
This shows that one may reduce the relation of isometry between Polish metric spaces to that of isometry between 
\textit{concave} Polish metric spaces.

\end{section}

%%%%%%%%%%%%%%%%%%%%%%%%%%%%%%%%%%%%%%%%%%%%%%%%%%%%%%%%%%%%%%%%%%%%%%%%%%%%%%%%%%
%%%%%%%%%%%%%%%%%%%%%%%%%%%%%%%%%%%%%%%%%%%%%%%%%%%%%%%%%%%%%%%%%%%%%%%%%%%%%%%%%%
\begin{section}{The computation} \label{complex}
As we saw above, 
the results of Weaver seem to imply that one may reduce the relation of dilatation between concave Polish metric spaces to 
the relation of isometry between separable Banach spaces. Furthermore, the relation of isometry 
between Polish spaces  reduces to that of isometry between concave metric spaces. \\
This gives us a good idea as to how one may Borel reduce the relation of isometry between Polish metric 
spaces to that of linear isometry 
between Banach separable spaces: one simply has to find a way to replace the word "dilatation" by isometry in the first sentence above 
(at least for a big enough class of Polish metric spaces), and then to show that all the reductions involved are Borel.\\

First, we need a way to code Polish metric spaces. Given what we saw in the introduction, 
it would be natural to consider ${\mathcal F}(\U)$, 
endowed with the Effros Borel structure, as being the (Borel) space of all Polish metric spaces. This is what Gao and Kechris did 
in \cite{gaokec}, where they showed that the relation $\backsimeq_i^{\U}$ induced by the action of $Iso(\U)$ on 
${\mathcal F}(\U)$ is universal for relations induced by a Borel action of some Polish group.

We shall indeed make use of this coding; to simplify the proof below, 
we introduce another way of coding Polish metric spaces.\\ 

We follow Vershik \cite{Vershik} and define the set of \textit{codes of Polish metric spaces} as the set ${\mathcal M}$ of all 
$\mathbf{d}=(d_{i,j}) \in \R^{\omega \times \omega}$ such that: \\
\begin{align*}
    (1) & \ \  \forall i,j  \ d_{i,j} \geq 0;  \\
    (2) & \ \  \forall i  \ d_{i,i}=0;  \\
    (3) & \ \  \forall i,j  \ d_{i,j}=d_{j,i};  \\
    (4) & \ \  \forall i,j,k  \ d_{i,k}\leq d_{i,j}+d_{j,k}.
\end{align*}

Then, given a code $\mathbf{d}=(d_{i,j})$, one may naturally define a distance on the quotient of $\N$ by $\sim$, 
where $(i \sim j) \Leftrightarrow (d_{i,j}=0)$; then we associate to ${\mathbf d}$ 
the Polish metric space $(P,d)$ which is the completion of 
$\N_{/_{\sim}}$ endowed with the distance induced by ${\mathbf d}$.

Notice that ${\mathcal M}$ is a closed subset of $\R^{\omega
\times \omega}$ (endowed with the product topology); thus, the induced topology turns it into a 
Polish metric space.\\
It is then easy to check (see \cite{clemens}) that the relation of isometry in the codes, which we denote by
 $\backsimeq_i^{c}$, is analytic. \\
Also, it is important to point out here that results of Clemens, Gao and Kechris (see \cite{clemensgaokechris}) imply that
$$ (\backsimeq_i^{c}) \sim_B (\backsimeq_i) \sim_B (\backsimeq_i^{\U}) \ \ .$$
This means that the two codings of Polish spaces that we have introduced are equivalent from our point of view. \\
To see this, notice that the theorem of Kuratowski-Ryll-Nardzewski enables one to show 
easily that $(\backsimeq_i) \leq_B (\backsimeq_i^c)$. \\
Also, since $\U$ is finitely injective, one may find a Borel map 
$\Theta$ which sends any code $\mathbf{d}$ onto a closed set $\Theta(\mathbf{d}) \subset \U$ 
which is isometric to the space coded by ${\mathbf d}$; this proves that $(\backsimeq_i^c) \leq_B (\backsimeq_i)$. \\
Furthermore, results of Gao and Kechris \cite{gaokec} show that we may also ensure that
$ \mathbf{d} \backsimeq_i^c \mathbf{d'}  \Leftrightarrow
 \Theta(\mathbf{d}) \backsimeq^{\U}_i \Theta(\mathbf{d'}) \ ,$ and that $(\backsimeq_i) \sim_B (\backsimeq_i^{\U})$.\\
Depending on the situation, one of these relations is more convenient to use than the others, 
which explains why we introduced the set ${\mathcal M}$ and the relation $\backsimeq_i^c$.\\

We now show how one may replace the word "dilatation" by isometry in the reasoning described 
in the beginning of this section.\\ 

If $(X,d_X)$ and $(Y,d_Y)$ are two metric spaces, then $(X,d_X)$ and $(Y,d_Y)$ are isometric if, and only if, 
$(X,\frac{d_X}{1+d_X})$ and $(Y,\frac{d_Y}{1+d_Y})$ are isometric. \\
Furthermore, the mapping "$d \mapsto \frac{d}{1+d}$" is continuous in the codes: hence, isometry between Polish metric spaces 
Borel reduces to isometry between bounded Polish metric spaces, and isometry between unbounded Polish metric spaces 
Borel reduces to isometry between Polish spaces of diameter (exactly) $1$.\\
It is also easy to associate continuously (in the codes) a code for an unbounded Polish metric space $\Psi({\mathbf d})$ 
to any code ${\mathbf d}$ for a bounded Polish metric space $X$, in such a way that 
$$(\Psi({\mathbf d}) \backsimeq_i \Psi({\mathbf d'})) \Leftrightarrow ({\mathbf d} \backsimeq_i^c {\mathbf d'}) \ .$$ 
Hence, the relation of isometry between Polish metric spaces Borel reduces to the relation of isometry between Polish spaces of 
diameter (exactly) $1$, so these relations are actually Borel bireducible.\\

The point of fixing the diameter of the spaces we consider is of course that, if $X$ and $Y$ are two bounded Polish metric spaces 
of the same diameter, then any dilatation of $X$ onto $Y$ actually must be an isometry.

Given the results of Gao and Kechris quoted above, we may find a Borel map $\Phi_0$,
defined on the sets of codes of Polish metric spaces of diameter (exactly) $1$, and taking 
its values in ${\mathcal F}(\U)$ (endowed with the Effros Borel structure), such that :\\

- If $\mathbf{d}$ codes a Polish metric space $(P,d)$, then
$\Phi_0(\mathbf{d})$ contains $0$, and is isometric to $(P,\sqrt{d})$. \\
- If $\mathbf{d} \backsimeq_i^c \mathbf{d'} $ then $\Phi_0(\mathbf{d}) \backsimeq_i^{\U} \Phi_0(\mathbf{d'})$.\\

Let $\Phi(\mathbf{d})$ denote the closed linear span of 
$\Phi_0(\mathbf{d})$ in $F(\U)$. Then we see that $\Phi \colon {\mathcal M} \to {\mathcal F}(F(\U))$ is Borel. \\

Then, thanks to the results of Weaver detailed above, we obtain the following result:

\begin{thm} Let $\mathbf{d}$, $\mathbf{d'}$ code two Polish metric spaces $P$, $P'$ of diameter (exactly) $1$. \\
Then the following assertions are equivalent:\\\
(1) $P$ and $P'$ are isometric. \\
(2) $\Phi(\mathbf{d})$ and $\Phi(\mathbf{d'})$ are isometric . \\
(3) There is a linear isometry of $F(\U)$ which maps
$\Phi(\mathbf{d})$ onto $\Phi(\mathbf{d'})$.\\
\end{thm}

{\bf Proof: } \\
$(2) \Rightarrow (1)$ is a direct consequence of the results of Weaver. 
Indeed, if $\Phi(\mathbf{d})$ and $\Phi(\mathbf{d})$ are isometric, then $\Phi_0(\mathbf{d})$
and $\Phi_0(\mathbf{d'})$ are isometric (since they are concave, and have the same diameter); then, 
by definition of $\Phi_0$, we know that $(P,\sqrt{d})$ and $(P',\sqrt{d'})$ are isometric, 
hence $P$ and $P'$ must be isometric too.\\

$(3) \Rightarrow(2)$ is a triviality, and $(1) \Rightarrow(3)$ is a consequence of the fact that any 
isometry of $\U$ mapping $\Phi_0(\mathbf{d})$
onto $\Phi_0(\mathbf{d'})$ extends to a linear isometry
of $F(\U)$ mapping the closed linear span of $\Phi_0(\mathbf{d})$
onto that of $\Phi_0(\mathbf{d'})$. $\hfill \lozenge$\\

This shows that $\Phi$ is a Borel reduction of isometry between Polish spaces of diameter (exactly) $1$ to the relation 
of linear isometry between Banach spaces. Given the results of Gao and Kechris \cite{gaokec} detailed in the introduction, 
this is enough to compute the exact complexity of the relation of isometry between separable Banach spaces.

\begin{thm}
The relation of isometry between separable Banach spaces is Borel bi-reducible to the universal relation for 
relations induced by a Borel action of a Polish group.
\end{thm}

Unfortunately, the results above do not seem to shed any light 
on the many questions remaining open in related areas; 
for instance we still do not know the exact complexity of isomorphism between separable 
Banach spaces, and have very little information about the relation of homeomorphism between compact metric spaces.\\
Also, it is worth pointing out that $F(\U)$ seems to be quite a remarkable object, as shown by the theorem of Holmes stated 
earlier. To our knowledge, next to nothing is known about the geometry of this space; it 
might prove fruitful to delve further in that direction (part of the problem here is that it is hard to ask "good" 
questions about the geometry of this space, since its definition makes the study of the "usual" questions quite 
difficult, at least for us).

\end{section}

\end{document}